\expandafter\ifx\csname mthreemacsloaded\endcsname\relax\else \fi

\magnification1100
\input amstex


 \catcode`\@=11
 \let\wlog@ld\wlog
 \def\wlog#1{\relax}

 \newif\ifIN@
 \def\m@rker{\m@@rker}
 \def\IN@{\expandafter\INN@\expandafter}
 \long\def\INN@0#1@#2@{\long\def\NI@##1#1##2##3\ENDNI@
    {\ifx\m@rker##2\IN@false\else\IN@true\fi}%
     \expandafter\NI@#2@@#1\m@rker\ENDNI@}
  \newtoks\Initialtoks@  \newtoks\Terminaltoks@
  \def\SPLIT@{\expandafter\SPLITT@\expandafter}
  \def\SPLITT@0#1@#2@{\def\TTILPS@##1#1##2@{%
     \Initialtoks@{##1}\Terminaltoks@{##2}}\expandafter\TTILPS@#2@}
  \newtoks\Trimtoks@

 \def\ForeTrim@{\expandafter\ForeTrim@@\expandafter}
 \def\ForePrim@0 #1@{\Trimtoks@{#1}}
 \def\ForeTrim@@0#1@{\IN@0\m@rker. @\m@rker.#1@%
     \ifIN@\ForePrim@0#1@%
     \else\Trimtoks@\expandafter{#1}\fi}
 
  \def\Trim@0#1@{%
      \ForeTrim@0#1@%
      \IN@0 @\the\Trimtoks@ @%
        \ifIN@
             \SPLIT@0 @\the\Trimtoks@ @\Trimtoks@\Initialtoks@
             \IN@0\the\Terminaltoks@ @ @%
                 \ifIN@
                 \else \Trimtoks@ {FigNameWithSpace}%
                 \fi
        \fi
      }

  \font\titlebold=cmbx12 scaled 1200
  \font\twelvebold=cmbx12
  \font\tenbold=cmbx10
  \font\ninebold=cmbx9
  \font\sevenbold=cmbx7
  \font\fivebold=cmbx5

  \input amssym.def \input amssym
     \font\titlemsa=msam10 at 14.4pt
     \font\titlemsb=msbm10 at 14.4pt
     \font\titleeufm=eufm10 at 14.4pt
     \font\twelvemsa=msam10 scaled 1200
     \font\twelvemsb=msbm10 scaled 1200
     \font\twelveeufm=eufm10 scaled 1200
     \font\ninemsa=msam9
     \font\ninemsb=msbm9
     \font\nineeufm=eufm9

   \ifx\cyrfam\undefined
   \else
     \immediate\write16{}%
     \message{ !!! cyr fonts already defined. !!! }
     \message{ --- edit out superfluous font defs? }
   \fi
   \newfam\cyrfam
       \font\titlecyr=wncyr10 scaled 1440 
       \font\twelvecyr=wncyr10 scaled 1200
       \font\tencyr=wncyr10
       \font\ninecyr=wncyr9
       \font\sevencyr=wncyr7
       \font\sixcyr=wncyr6

   \newfam\eusmfam
       \font\titleeusm=eusm10 scaled 1440
       \font\twelveeusm=eusm10 scaled 1200
       \font\teneusm=eusm10
       \font\nineeusm=eusm9
       \font\seveneusm=eusm7
       
       \font\fiveeusm=eusm5

\let\Cal\cal

    \font\ninemrm=cmr9 
    \font\ninei=cmmi9
    \font\ninesy=cmsy9 
    \skewchar\ninei='177
    \skewchar\ninesy='60

  \font\twelvemrm=cmr10 at 12pt 
  \font\twelvei=cmmi10 at 12pt
  \font\twelvesy=cmsy10 at 12pt

  \font\titlemrm=cmr10 at 14.4pt 
  \font\titlei=cmmi10 at 14.4pt
  \font\titlesy=cmsy10 at 14.4pt


  \def\Smallfonts{\ninepoint}

  \def\Hfont{\titlepoint\bf}
  \def\Authorfont{\twelvepoint\it}
  \def\HHfont{\twelvepoint\bf}
  \def\HHHfont{\bf}
  \def\Bibfont{\tenbf}
  \def\Coordfont{\nineit }

  \def \thfont {\bf }
  \def \pffont {\it\itSpacing }
  \def \rkfont {\bf }
  \def \dffont {\bf }
  \def \egfont {\bf }

 \def\ninepoint{%
  \def\rm{\fam0\ninerm}%
    \textfont0=\ninemrm  \scriptfont0=\sevenrm  \scriptscriptfont0=\fiverm
    \textfont1=\ninei    \scriptfont1=\seveni   \scriptscriptfont1=\fivei
  \def\mit{\fam1\ninei}%
  \def\oldstyle{\fam1\ninei}%
    \textfont2=\ninesy   \scriptfont2=\sevensy  \scriptscriptfont2=\fivesy
    \textfont3=\tenex    \scriptfont3=\tenex    \scriptscriptfont3=\tenex
  \def\it{\fam\itfam\nineit}%
    \textfont\itfam=\nineit
  \def\bf{\ifmmode\fam\bffam\else\ninebf\fi}%
    \textfont\bffam=\ninebold 
    \scriptfont\bffam=\sevenbold 
    \scriptscriptfont\bffam=\fivebold%
  \def\msa{\fam\msafam\ninemsa}%
    \textfont\msafam=\ninemsa 
    \scriptfont\msafam=\sevenmsa
    \scriptscriptfont\msafam=\fivemsa%
  \def\msb{\fam\msbfam\ninemsb}%
    \textfont\msbfam=\ninemsb%
    \scriptfont\msbfam=\sevenmsb%
    \scriptscriptfont\msbfam=\fivemsb%
  \def\eufm{\fam\eufmfam\nineeufm}%
    \textfont\eufmfam=\nineeufm
    \scriptfont\eufmfam=\seveneufm
    \scriptscriptfont\eufmfam=\fiveeufm
   \def\eusm{\fam\eusmfam\nineeusm}%
     \textfont\eusmfam=\nineeusm
     \scriptfont\eusmfam=\seveneusm
     \scriptscriptfont\eusmfam=\fiveeusm
   \def\cyr{\fam\cyrfam\ninecyr}%
     \textfont\cyrfam=\ninecyr
     \scriptfont\cyrfam=\sevencyr
     \scriptscriptfont\cyrfam=\sixcyr
  \setbox\strutbox=\hbox{\vrule
      height7pt depth3pt width0pt}%
   \baselineskip=10.8pt\rm}

 \let\eightpoint\ninepoint 

 \def\tenpoint{%
  \def\rm{\fam0\tenrm}%
    \textfont0=\tenmrm \scriptfont0=\sevenrm \scriptscriptfont0=\fiverm%
  \def\mit{\fam1\teni}%
  \def\oldstyle{\fam1\teni}%
    \textfont1=\teni   \scriptfont1=\seveni  \scriptscriptfont1=\fivei%
    \textfont2=\tensy  \scriptfont2=\sevensy \scriptscriptfont2=\fivesy%
    \textfont3=\tenex  \scriptfont3=\tenex   \scriptscriptfont3=\tenex%
  \def\it{\fam\itfam\tenit}%
    \textfont\itfam=\tenit%
  \def\bf{\ifmmode\fam\bffam\else\tenbf\fi}%
    \textfont\bffam=\tenbold
    \scriptfont\bffam=\sevenbold%
    \scriptscriptfont\bffam=\fivebold%
  \def\msa{\fam\msafam\tenmsa}%
    \textfont\msafam=\tenmsa%
    \scriptfont\msafam=\sevenmsa%
    \scriptscriptfont\msafam=\fivemsa%
  \def\msb{\fam\msbfam\tenmsb}%
    \textfont\msbfam=\tenmsb%
    \scriptfont\msbfam=\sevenmsb%
    \scriptscriptfont\msbfam=\fivemsb%
  \def\eufm{\fam\eufmfam\teneufm}%
   \textfont\eufmfam=\teneufm
   \scriptfont\eufmfam=\seveneufm
   \scriptscriptfont\eufmfam=\fiveeufm
   \def\eusm{\fam\eusmfam\teneusm}%
    \textfont\eusmfam=\teneusm
    \scriptfont\eusmfam=\seveneusm
    \scriptscriptfont\eusmfam=\fiveeusm
   \def\cyr{\fam\cyrfam\tencyr}%
    \textfont\cyrfam=\tencyr
    \scriptfont\cyrfam=\sevencyr
    \scriptscriptfont\cyrfam=\sixcyr
  \setbox\strutbox=\hbox{\vrule %
      height8.5pt depth3.5ptwidth0pt}%
  \baselineskip=\StdBaselineskip\rm}

 \def\twelvepoint{%
  \def\rm{\fam0\twelverm}%
    \textfont0=\twelvemrm \scriptfont0=\tenmrm \scriptscriptfont0=\sevenrm
    \textfont1=\twelvei   \scriptfont1=\teni   \scriptscriptfont1=\seveni
  \def\mit{\fam1\twelvei}%
  \def\oldstyle{\fam1\twelvei}%
    \textfont2=\twelvesy  \scriptfont2=\tensy  \scriptscriptfont2=\sevensy
    \textfont3=\tenex  \scriptfont3=\tenex  \scriptscriptfont3=\tenex
  \def\it{\fam\itfam\twelveit}%
    \textfont\itfam=\twelveit
  \def\bf{\ifmmode\fam\bffam\else\twelvebf\fi}%
    \textfont\bffam=\twelvebold
    \scriptfont\bffam=\tenbold%
    \scriptscriptfont\bffam=\sevenbold%
  \def\msa{\fam\msafam\twelvemsa}%
    \textfont\msafam=\twelvemsa%
    \scriptfont\msafam=\tenmsa%
    \scriptscriptfont\msafam=\sevenmsa%
  \def\msb{\fam\msbfam\twelvemsb}%
    \textfont\msbfam=\twelvemsb%
    \scriptfont\msbfam=\tenmsb%
    \scriptscriptfont\msbfam=\sevenmsb%
  \def\eufm{\fam\eufmfam\twelveeufm}%
   \textfont\eufmfam=\twelveeufm
   \scriptfont\eufmfam=\teneufm
   \scriptscriptfont\eufmfam=\seveneufm
   \def\eusm{\fam\eusmfam\twelveeusm}%
    \textfont\eusmfam=\twelveeusm
    \scriptfont\eusmfam=\teneusm
    \scriptscriptfont\eusmfam=\seveneusm
   \def\cyr{\fam\cyrfam\tencyr}%
    \textfont\cyrfam=\twelvecyr
    \scriptfont\cyrfam=\tencyr
    \scriptscriptfont\cyrfam=\sevencyr
  \setbox\strutbox=\hbox{\vrule
      height10.2pt depth4.55pt width0pt}%
  \baselineskip=14pt\rm}

 \def\titlepoint{%
    \textfont0=\titlemrm \scriptfont0=\twelvemrm \scriptscriptfont0=\tenmrm
    \textfont1=\titlei   \scriptfont1=\twelvei   \scriptscriptfont1=\teni
  \def\mit{\fam1\titlei}%
  \def\oldstyle{\fam1\titlei}%
    \textfont2=\titlesy  \scriptfont2=\twelvesy  \scriptscriptfont2=\tensy
    \textfont3=\tenex
    \scriptfont3=\tenex
    \scriptscriptfont3=\tenex
  \def\it{\fam\itfam\titleit}%
    \textfont\itfam=\titleit
  \def\bf{\ifmmode\fam\bffam\else\titlebf\fi}%
    \textfont\bffam=\titlebold
    \scriptfont\bffam=\twelvebold%
    \scriptscriptfont\bffam=\tenbold%
  \def\msa{\fam\msafam\titlemsa}%
    \textfont\msafam=\titlemsa%
    \scriptfont\msafam=\twelvemsa%
    \scriptscriptfont\msafam=\tenmsa%
  \def\msb{\fam\msbfam\titlemsb}%
    \textfont\msbfam=\titlemsb%
    \scriptfont\msbfam=\twelvemsb%
    \scriptscriptfont\msbfam=\tenmsb%
  \def\eufm{\fam\eufmfam\titleeufm}%
    \textfont\eufmfam=\titleeufm
    \scriptfont\eufmfam=\twelveeufm
    \scriptscriptfont\eufmfam=\teneufm
   \def\eusm{\fam\eusmfam\titleeusm}%
     \textfont\eusmfam=\titleeusm
     \scriptfont\eusmfam=\twelveeusm
     \scriptscriptfont\eusmfam=\teneusm
   \def\cyr{\fam\cyrfam\tencyr}%
    \textfont\cyrfam=\titlecyr
    \scriptfont\cyrfam=\twelvecyr
    \scriptscriptfont\cyrfam=\tencyr
  \setbox\strutbox=\hbox{\vrule
      height12.3pt depth5.54pt width0pt}%
  \baselineskip=16pt\rm}

\newbox\AuthorBox\newbox\TitleBox
\newbox\TFLinebox
\newbox\FLinebox
\newbox\HLinebox
\def\SetTFLinebox#1{\setbox\TFLinebox=\hbox{#1}}
\def\SetFLinebox#1{\setbox\FLinebox=\hbox{#1}}
\def\SetHLinebox#1{\setbox\HLinebox=\hbox{#1}}

 \def\SetAuthorHead#1{%
     \setbox\AuthorBox=\hbox{\ninepoint \it 
           \ignorespaces\frenchspacing#1\unskip}}
 \def\SetTitleHead#1{%
     \setbox\TitleBox=\hbox{\ninepoint \it
           \ignorespaces\frenchspacing#1\unskip}}

  \def\itSpacing{\relax}
  \def\itSpacingOff{\relax}


 \def\Hrule{\hrule width0pt height0pt}

  \newskip\ProcSkip \ProcSkip 8pt plus2pt minus2pt

 \newskip\LastSkip
 \def\SaveLastSkip{\LastSkip\lastskip}
 \def\RestoreLastSkip{\vskip-\LastSkip\vskip\LastSkip}

 \def\NoindentAfter{\everypar={\setbox0=\lastbox\everypar={}}}

 \long\def\H#1\par#2\par{\notenumber=0 \titlepagetrue%
    {
    \baselineskip=20pt
    \parindent=0pt\parskip=0pt\frenchspacing
    \leftskip=0pt plus .2\hsize minus .3\hsize
    \rightskip=0pt plus .2\hsize minus .3\hsize
 \def\\{\unskip\break}%
    \pretolerance=10000 \Hfont #1\unskip\break
     \vskip7pt\Hrule
\hfill \Authorfont #2\hfill\hfill\unskip}
    \vskip48pt plus 4pt minus 4pt
    \par\NoindentAfter\rm}

 \long\def\Hi#1\par#2\par{\notenumber=0 \titlepagetrue%
    {  \baselineskip=0pt  \parindent=0pt\parskip=0pt\frenchspacing
    \leftskip=0pt plus .2\hsize minus .3\hsize
    \rightskip=0pt plus .2\hsize minus .3\hsize
}
    \rm}


 \newdimen\PageRemainder
  \def\SetPageRemainder{
     \PageRemainder=\pagegoal
     \ifdim\PageRemainder=\maxdimen\PageRemainder=\vsize
     \else\advance\PageRemainder by -1\pagetotal\fi}

  \def\Rpt@{}\def\Rpt@@{}

  \long\def\HH#1\par{\par
  \SaveLastSkip\removelastskip\goodbreak
  \ifdim\LastSkip<30pt 
     \LastSkip 30pt
plus 3pt minus 2pt\fi
  \SetPageRemainder\advance\PageRemainder-\LastSkip
  \ifdim\PageRemainder<150pt
       \edef\Rpt@{remain = \the\PageRemainder\noexpand\\
                pagetotal=\the\pagetotal\noexpand\\
                           pagegoal=\the\pagegoal}%
          \fi
   \ifdim\PageRemainder<65pt 
       \ifdim\PageRemainder > 0pt
          \edef\Rpt@@{\noexpand\\
                      Had HH PageRemainder$<$\relax 65pt\noexpand\\
                      Hence forced break!}%
     \vskip 0pt plus .2\PageRemainder\eject 
    \fi\fi
    \vskip\LastSkip\Hrule 
    \pretolerance=10000\rightskip=0pt plus 3em
    \hangafter1 \hangindent=2.2em%
    \noindent
    \HHfont \unskip \Ednote{\Rpt@\Rpt@@}%
            \def\Rpt@{}\def\Rpt@@{}%
            \ignorespaces
            #1\par\rightskip=0pt\pretolerance=\StdPretolerance%
    \NoindentAfter
\tenpoint\rm%
     \medskip \vskip\ProcSkip}

  \long\def\HHH#1\par{\par%
  \SaveLastSkip\removelastskip\goodbreak
  \ifdim\LastSkip<\ProcSkip%
     \LastSkip\ProcSkip\fi
  \SetPageRemainder\advance\PageRemainder-\LastSkip
  \ifdim\PageRemainder<150pt
       \edef\Rpt@{remain = \the\PageRemainder\noexpand\\
                pagetotal=\the\pagetotal\noexpand\\
                           pagegoal=\the\pagegoal}%
       \fi
   \ifdim\PageRemainder<48pt  
        \ifdim\PageRemainder > 0pt
             \edef\Rpt@@{\noexpand\\
                      Had HHH PageRemainder$<$\relax48pt\noexpand\\
                      Hence forced break!}%
       \vskip 0pt plus .2\PageRemainder\eject 
      \fi\fi
   \vskip\LastSkip\par\noindent
   \HHHfont \unskip\Ednote{\Rpt@\Rpt@@}%
  \def\Rpt@{}\def\Rpt@@{}%
  \ignorespaces
   #1\unskip.\quad\rm\ignorespaces
   \ignorepars}

  \long\def\ignorepars#1\par{\def\Test{#1}%
     \ifx\Test\Empty\def\This{\ignorepars}%
        \else\def\This{\Test\par}\fi
           \This}
  \def\Empty{}

 \def\Abstract#1\par{\bgroup\Smallfonts\narrower\HHH #1\par}
 \def\endAbstract{\par\egroup}


 \def\ProcBreak{\par%
    \ifdim\lastskip<8pt%
    \removelastskip%
    \penalty-200\vskip\ProcSkip\fi}

 \def\th#1\par{\ProcBreak \noindent
   {\thfont\ignorespaces
    #1\unskip.}\it\itSpacing\kern.4em\ignorepars}

 \def\endth{\ProcBreak\rm\itSpacingOff }


 \def\pf#1\par{\ProcBreak %
    \noindent\pffont#1\unskip.\rm\itSpacingOff{\kern .7em}\ignorepars}

 \def\endpf{\medskip \ProcBreak } 

  \def\qedbox{\hbox{\vbox{
    \hrule width0.2cm height0.2pt
    \hbox to 0.2cm{\vrule height 0.2cm width 0.2pt
             \hfil\vrule height0.2cm width 0.2pt}
    \hrule width0.2cm height 0.2pt}\kern1pt}}

  \def\qed{\ifmmode\qedbox
    \else\unskip\ \hglue0mm\hfill\qedbox\ProcBreak\fi}

  \def \rk #1\par{\ProcBreak
     \noindent{\rkfont\ignorespaces #1\unskip.}%
     \rm\kern.6em\ignorepars}

  \def \endrk {\medskip\ProcBreak }

  \def \df #1\par{\ProcBreak
     \noindent{\dffont\unskip\ignorespaces #1\unskip.}%
     \rm\kern.6em\ignorepars}

  \def \enddf {\medskip\ProcBreak }

  \def \eg #1\par{\ProcBreak
     \noindent\egfont\unskip\ignorespaces #1\unskip.
     \rm\kern.6em\ignorepars}

  \newdimen\Overhang

   \def\MaxTag@#1#2#3#4#5{\setbox0=\hbox{#4\ignorespaces#2\unskip}%
     \dimen0=\wd0\advance\dimen0 by#3
     \ifdim\dimen0<#5\relax\dimen0=#5\fi
     \expandafter\edef\csname #1Hang\endcsname{\the\dimen0}}

 \def\MaxItemTag#1{\MaxTag@{Item}{#1}{.4em}{\ItemStyle}{\parindent}}%
 \def\MaxItemItemTag#1{%
        \MaxTag@{ItemItem}{#1}{.4em}{\ItemItemStyle}{\parindent}}
 \def\MaxNrTag#1{\MaxTag@{Nr}{#1}{.5em}{\NrStyle}{\parindent}}
 \def\MaxReferenceTag#1{%
        \MaxTag@{Reference}{[#1]}{.6em}{\ninerm}{\parindent}}
 \def\MaxFootTag#1{\MaxTag@{Foot}{#1}{.4em}{\ninerm}{\z@}}

  \def\SetOverhang@{\Overhang=.8\dimen0%
     \advance\Overhang by \wd0\relax
     \ifdim\Overhang>\hangindent\relax
       \advance\Overhang by .25\dimen0%
       \Ednote{Tag is pushing text.}\osumess{Tag is pushing text.}%
     \else\Overhang=\hangindent
     \fi}

   \def\Item#1{\par\noindent
      \hangafter1\hangindent=\ItemHang
      \setbox0=\hbox{\ItemStyle\ignorespaces#1\unskip}%
      \dimen0=.4em\SetOverhang@
      \rlap{\box0}\kern\Overhang\ignorespaces}

   \def\ItemItem#1{\par\noindent
      \hangafter1\hangindent=\ItemItemHang
      \setbox0=\hbox{\ItemItemStyle\ignorespaces#1\unskip}%
      \dimen0=.4em\SetOverhang@
      \advance\hangindent by \ItemHang
      \kern\ItemHang\rlap{\box0}%
      \kern\Overhang\ignorespaces}

  \def\Nr#1{\par\noindent\hangindent=\NrHang 
    \setbox0=\hbox{\NrStyle\ignorespaces#1\unskip}%
    \dimen0=.5em\SetOverhang@
    \rlap{\box0}\kern\Overhang
    \hangindent=\z@\ignorespaces}

   \newskip\Rosterskip\Rosterskip 1pt plus1pt 
   \def\Roster{\par\ifdim\lastskip<\Rosterskip\removelastskip\vskip\Rosterskip\fi
    \bgroup}
   \def\endRoster{\par\global\edef\LastSkip@{\the\lastskip}\removelastskip
       \egroup\penalty-50\LastSkip\LastSkip@\relax
       \ifdim\LastSkip<\Rosterskip\LastSkip\Rosterskip\fi
       \vskip\LastSkip}




 \def\cite#1{
    \def\nextiii@##1,##2\end@{{\frenchspacing\rm 
      \lBr\ignorespaces##1\unskip{\rm,~\ignorespaces##2}\rBr}}%
    \IN@0,@#1@%
    \ifIN@\def\next{\nextiii@#1\end@}\else
    \def\next{{\rm\lBr#1\rBr}}\fi\next}


   \def \Bib#1\par{%
       \par\removelastskip\SetPageRemainder
       \ifdim\PageRemainder < 97pt
        \ifdim\PageRemainder > 0pt
        \vfill\eject
       \fi\fi
    \ProcBreak \par\begingroup\parskip=0 pt%
    \goodbreak \vskip 15 pt plus 10 pt
    \noindent\null\hfill\Bibfont
      \ignorespaces #1\unskip\hfill\null\par 
    \frenchspacing \Smallfonts\rm
    \parskip=2.5 pt plus 1 pt minus.5pt%
    \nobreak\vskip 12pt plus 2pt minus2pt\nobreak
    \leftskip=0 pt \baselineskip=10.5pt}

 \def\ReferenceTagSlide{0em}
  \def\ReferenceTagGap{.5em}

  \def \rf#1{\par\noindent
     \hangafter1\hangindent=\ReferenceHang      
     \setbox0=\hbox{\ninerm[\ignorespaces#1\unskip]}%
     \dimen0=\ReferenceTagGap\SetOverhang@
     \rlap{\kern\ReferenceTagSlide\box0}%
     \kern\Overhang\ignorespaces}

  \def\ref#1\par#2\par#3\par#4\par{%
     \rf{#1}#2\unskip,\ #3\unskip,\
     #4\unskip.}

  \def\endBib{\par\endgroup\vskip 12pt minus 6pt }


  \long\def\Coordinates#1\endCoordinates{
 {\par\vskip4pt\def\\{\unskip, }\Coordfont\baselineskip10.5pt\noindent#1}}

 \def\pagecontents{
  \gdef\Pagetot@l{\pagetotal}
  \ifvoid\TRMargIns\else
    \rlap{\kern\hsize\kern10pt\vbox to 0pt{%
         \box\TRMargIns\vss}}\fi
  \ifvoid\topins\else\unvbox\topins\fi
   \dimen@=\dp\@cclv \unvbox\@cclv 
   \ifvoid\footins\else 
     \vskip\skip\footins
     \footnoterule
     \unvbox\footins\fi
   \ifr@ggedbottom \kern-\dimen@ \vfil \fi}


 \newcount\Ht 

 \def \Acc{\expandafter } 

 \def\swthat{\raise -1.1 ex\hbox{\sam$\widehat{}$}}
 \def\swttilde{\raise -1.2 ex\hbox{\sam$\widetilde{}$}}
 \def \overdot{{\raise .2 ex \hbox to 0pt {\hss\bf\smash{.}\hss}}}
 \def \overcircle{{\raise .1 ex \hbox to 0pt
    {\sam$\eightpoint\scriptstyle\hss\circ\hss$}}}

 \def \Mathaccent#1#2{{\sam 
  \setbox4=\hbox{$\vphantom{#2}$}
  \Ht=\ht4 
  \setbox5=\hbox{${#1}$}
  \setbox6=\hbox{${#2}$}
  \setbox7=\hbox to .5\wd6{}
  \copy7\kern .1\Ht \raise\Ht sp\hbox{\copy5}\kern-.1\Ht
  \copy7\llap{\box6}
  }}

  \def\SwtCheck #1{
        \ifmmode \check{#1}%
                \else \v {#1}%
                \fi}

 \def\barpartial {%
   \kern .17 em
    \overline {\kern -.17 em\partial\kern-.03 em}%
    \kern .03 em}

 
  \def\Overline#1{\setbox1=\hbox{\sam ${#1}$}%
      \ifdim \wd1 > 6pt
    \kern .11 em
    \overline {\kern -.11 em#1\kern-.14 em}
    \kern .14 em
  \else
    \kern .03 em
    \overline {\kern -.03 em#1\kern-.04 em}
    \kern .04 em
  \fi}

 \def\SOverline#1{\setbox1=\hbox{\sam ${#1}$}%
      \ifdim \wd1 > 7pt
    \kern .22 em
    \overline {\kern -.22 em#1\kern-.09 em}%
    \kern .09 em
  \else
    \kern .10 em
    \overline {\kern -.10 em#1\kern-.04 em}%
    \kern .04 em
  \fi}


 \def\Underline#1{\setbox1=\hbox{\sam ${#1}$}%
      \ifdim \wd1 > 6pt
    \kern .11 em
    \underline {\kern -.11 em#1\kern-.14 em}
    \kern .14 em
  \else
    \kern .03 em
    \underline {\kern -.03 em#1\kern-.04 em}
    \kern .04 em
  \fi}

 \def\SUnderline#1{\setbox1=\hbox{\sam ${#1}$}%
      \ifdim \wd1 > 7pt
    \kern .04 em
    \underline {\kern -.04 em#1\kern-.2 em}%
    \kern .2 em
  \else
    \kern .0 em
    \underline {\kern -.0 em#1\kern-.15 em}%
    \kern .15 em
  \fi}


 \def \Blackbox
   {\leavevmode\hskip .3pt \vbox
   {\hrule height 5pt\hbox{\hskip 4.5pt}}\hskip .5pt}

 \def \XX{\Blackbox\kern.5pt\Blackbox} 

  \def\.{.\kern1pt}

    \def\Hyphen{\edef\this{\the\hyphenchar\font}%
          \hyphenchar\font=-1\char\this\hyphenchar\font=\this}

 \ifx\undefined\text
  \def\text#1{\hbox{\rm #1}}\fi 



   \everymath{}  

  \def\PassMath@@{\aftergroup\AfterMath@} 

 \let\PassMath@\PassMath@@

 \def\AfterMath@{\futurelet\next\AfterMathMole@}

 \def\AfterMathMole@{
      \ifcat\next\space
          \def\this{}
      \else
      \ifcat\next\egroup %
        \def\this{\osumess{Handset mathsurround?? ---(see dollar brace)}}%
      \else
      \def\this{\AAfterMath@}
      \fi\fi
      \this}

 \def\hyphen@{-}
 \def\paren@{)}
 \def\apostr@{'}

 \def\MSC#1{\kern-.8\mathsurround#1\kern.8\mathsurround}

 \def\AAfterMath@#1{\def\Next{#1}
    \IN@0\Next @,.;:!?\relax @%
    \ifIN@\def\this{\MSC{\Next}}%
    \else
    \ifx\Next\hyphen@\def\this{\futurelet\next\AfterHyphen@}%
    \else
    \ifx\Next\paren@\def\this{#1}%
    \else 
    \ifx\Next\apostr@\def\this{#1}%
    \else \def\this{\osumess{Handset mathsurround??}%
                 #1}\fi\fi\fi\fi
    \this}

 \def\AfterHyphen@#1{\def\Next{#1}%
   \ifx\Next\hyphen@\def\this{--}\else
   \ifcat\next\space%
   \def\this{\kern-\mathsurround\kern.05em- \Next}\else
   \def\this{\kern-\mathsurround\kern.05em\Hyphen\Next}\fi\fi\this}

 \def\sam{\mathsurround=\z@\let\PassMath@\relax}  %
 \def\mas{\mathsurround=\StdMathsurround\let\PassMath@\PassMath@@}
 
 \def\Mas{\mathsurround=\StdMathsurround
                \everymath{\PassMath@}\let\PassMath@\PassMath@@}

 \def\m@th{\mathsurround=\z@\everymath{}}

 \def\m@@th{\mathsurround=\z@\everymath={}\let\m@th\relax}

\def\underbar#1{$\setbox\z@\hbox{#1}\dp\z@\z@
      \m@th \underline{\box\z@}$\relax}

\def\mathhexbox#1#2#3{\leavevmode
  \hbox{\m@@th$\m@th \mathchar"#1#2#3$}}

\def\dots{\relax\ifmmode\ldots\else$\m@th\ldots\,$\relax\fi}

\def\dotfill{\cleaders\hbox{\m@@th$\m@th \mkern1.5mu.\mkern1.5mu$}\hfill}
\def\rightarrowfill{$\m@th\mathord-\mkern-6mu%
  \cleaders\hbox{\m@@th$\mkern-2mu\mathord-\mkern-2mu$}\hfill
  \mkern-6mu\mathord\rightarrow$\relax}
\def\leftarrowfill{$\m@th\mathord\leftarrow\mkern-6mu%
  \cleaders\hbox{\m@@th$\mkern-2mu\mathord-\mkern-2mu$}\hfill
  \mkern-6mu\mathord-$\relax}

\def\downbracefill{$\m@th\braceld\leaders\vrule\hfill\braceru
  \bracelu\leaders\vrule\hfill\bracerd$\relax}
\def\upbracefill{$\m@th\bracelu\leaders\vrule\hfill\bracerd
  \braceld\leaders\vrule\hfill\braceru$\relax}

\def\angle{{\vbox{\m@@th\ialign{$\m@th\scriptstyle##$\crcr
      \not\mathrel{\mkern14mu}\crcr
      \noalign{\nointerlineskip}
      \mkern2.5mu\leaders\hrule height.34pt\hfill\mkern2.5mu\crcr}}}}

\def\big#1{{\m@@th\hbox{$\left#1\vbox to8.5\p@{}\right.\n@space$}}}
\def\Big#1{{\m@@th\hbox{$\left#1\vbox to11.5\p@{}\right.\n@space$}}}
\def\bigg#1{{\m@@th\hbox{$\left#1\vbox to14.5\p@{}\right.\n@space$}}}
\def\Bigg#1{{\m@@th\hbox{$\left#1\vbox to17.5\p@{}\right.\n@space$}}}
\def\n@space{\nulldelimiterspace\z@ \m@th}

\def\root#1\of{\setbox\rootbox\hbox{\m@@th$\m@th\scriptscriptstyle{#1}$}
  \mathpalette\r@@t}
\def\r@@t#1#2{\setbox\z@\hbox{\m@@th$\m@th#1\sqrt{#2}$\relax}
  \dimen@\ht\z@ \advance\dimen@-\dp\z@
  \mkern5mu\raise.6\dimen@\copy\rootbox \mkern-10mu \box\z@}

\def\mathph@nt#1#2{\setbox\z@\hbox{\m@@th$\m@th#1{#2}$}\finph@nt}

\def\mathsm@sh#1#2{\setbox\z@\hbox{\m@@th$\m@th#1{#2}$}\finsm@sh}

\def\@vereq#1#2{\lower.5\p@\vbox{\m@@th\baselineskip\z@skip\lineskip-.5\p@
    \ialign{$\m@th#1\hfil##\hfil$\crcr#2\crcr=\crcr}}}

\def\mathpalette#1#2{\sam\mathchoice{#1\displaystyle{#2}}%
  {#1\textstyle{#2}}{#1\scriptstyle{#2}}{#1\scriptscriptstyle{#2}}\mas}

\def\widehat#1{\setbox\z@\hbox{\sam$#1$}%
 \ifdim\wd\z@>\tw@ em\mathaccent"0\msbfam@5B{#1}%
 \else\mathaccent"0362{#1}\fi}
\def\widetilde#1{\setbox\z@\hbox{\sam$#1$}%
 \ifdim\wd\z@>\tw@ em\mathaccent"0\msbfam@5D{#1}%
 \else\mathaccent"0365{#1}\fi}

 \def\dots{\relax{}
  \ifmmode\def\thedots{\mdots@}\else\def\thedots{\tdots@}\fi %
  \thedots}

 \let\@ldeqno\eqno\let\@ldleqno\leqno
 \def\eqno{\everymath{}\@ldeqno} \def\leqno{\everymath{}\@ldleqno}

  \let\@ldeqalignno\eqalignno
  \def\eqalignno#1{\sam\@ldeqalignno{#1}\mas}
  \let\@ldeqalign\eqalign
  \def\eqalign#1{\sam\@ldeqalign{#1}\mas}

 \def\overrightarrow#1{\vbox{\m@th\ialign{##\crcr
      \rightarrowfill\crcr\noalign{\kern-\p@\nointerlineskip}
      $\hfil\displaystyle{#1}\hfil$\crcr}}}
 \def\overleftarrow#1{\vbox{\m@th\ialign{##\crcr
      \leftarrowfill\crcr\noalign{\kern-\p@\nointerlineskip}
      $\hfil\displaystyle{#1}\hfil$\crcr}}}
 \def\overbrace#1{\mathop{\vbox{\m@th\ialign{##\crcr\noalign{\kern3\p@}
      \downbracefill\crcr\noalign{\kern3\p@\nointerlineskip}
      $\hfil\displaystyle{#1}\hfil$\crcr}}}\limits}
 \def\underbrace#1{\mathop{\vtop{\m@th\ialign{##\crcr
      $\hfil\displaystyle{#1}\hfil$\crcr\noalign{\kern3\p@\nointerlineskip}
      \upbracefill\crcr\noalign{\kern3\p@}}}}\limits}

  \let\@ldmatrix\matrix
  \let\end@ldmatrix\endmatrix
  \def\matrix{\sam\@ldmatrix}
  \def\endmatrix{\end@ldmatrix\mas}
  \let\@ldgather\gather
  \let\end@ldgather\endgather
  \def\gather{\sam\@ldgather}
  \def\endgather{\end@ldgather\mas}
  \let\@ldalign\align
  \let\end@ldalign\endalign
  \def\align{\sam\@ldalign}
  \def\endalign{\end@ldalign\mas}
  \let\@ldaligned\aligned
  \let\end@ldaligned\endaligned
  \def\aligned{\sam\@ldaligned}
  \def\endaligned{\end@ldaligned\mas}
  \let\@ldtag\tag
  \def\tag{\sam\@ldtag}
   %

   \let\MinCDArrowWidth\minCDaw@




\newskip\insertskipamount\newskip\inserthardskipamount
\insertskipamount 6pt plus2pt 
\inserthardskipamount 6pt
\def\insertskip{\vskip\insertskipamount}
\newcount\SplitTest
\def\SetSplitTest{\SplitTest\insertpenalties
  \insert\topins{\floatingpenalty1}%
  \advance\SplitTest-\insertpenalties}
\def\midinsert{\par
 \SaveLastSkip\penalty-150\SetSplitTest\RestoreLastSkip
 \ifnum\SplitTest=-1
  \@midfalse\p@gefalse\else\@midtrue\fi\@ins}
\def\@ins{\par\begingroup\setbox\z@\vbox\bgroup%
  \vglue\inserthardskipamount}
\def\endinsert{\egroup 
  \if@mid \dimen@\ht\z@ \advance\dimen@\dp\z@
    \advance\dimen@\insertskipamount
    \advance\dimen@\pagetotal\advance\dimen@-\pageshrink
    \ifdim\dimen@>\pagegoal\@midfalse\p@gefalse\fi\fi
  \if@mid%
    \ifdim\lastskip<\insertskipamount\removelastskip\insertskip\fi
    \nointerlineskip\box\z@\penalty-200\insertskip
  \else%
    \SaveLastSkip
    \insert\topins{\penalty100 
    \splittopskip\z@skip
    \splitmaxdepth\maxdimen \floatingpenalty\z@
    \ifp@ge \dimen@\dp\z@
    \vbox to\vsize{\unvbox\z@\kern-\dimen@}
    \else \box\z@\nobreak\insertskip\fi}
    \RestoreLastSkip
   \fi\endgroup}


  \newcount\notenumber
  
  \def\note{\advance\notenumber by 1
    \footnote{\the\notenumber)}}

  \newbox\footbox

  \def\footnote#1{\let\@sf\empty
    \ifhmode\edef\@sf{\spacefactor\the\spacefactor}\/\fi
    \sam${}^{\fam0 #1}$\@sf\vfootnote{#1}}%

  \def\vfootnote#1{\insert\footins\bgroup
     \interlinepenalty100 \splittopskip=1pt
     \floatingpenalty=20000
     \leftskip=0pt\rightskip=0pt%
     \parindent=.3em
     \Smallfonts\rm
     \FootItem@{#1}
     \futurelet\next\fo@t}

  \def\FootItem@#1{\par\hangafter1\hangindent=\FootHang
     \setbox0=\hbox{\ignorespaces#1\unskip}%
     \dimen0=.4em\SetOverhang@
     \noindent\rlap{\box0}\kern\Overhang\ignorespaces}


  \def\fo@t{\ifcat\bgroup\noexpand\next \let\next\f@@t
    \else\let\next\f@t\fi \next}
  \def\f@@t{\bgroup\aftergroup\@foot\let\next}
  \def\f@t#1{\baselineskip=10pt\lineskip=1pt
            \lineskiplimit=0pt #1\@foot}%
  \def\@foot{
        \hbox{\vrule height0pt depth5pt width0pt}
        \egroup}
  \skip\footins=12 pt plus 0pt minus 0pt 
  \count\footins=1000 
  \dimen\footins=8in 



 \def\osumess#1{\EdSpider{\immediate\write16{Line \the\inputlineno: #1}}}%
 \def\HideEdStuff{\gdef\EdSpider##1{}}

 \font\BigSym=cmmi10 scaled \magstep 4

 \def\change{\InLMargin{\hbox{\BigSym \char63\kern10pt}}}

 \def\beginchange{\InLMargin{\hbox{\sam\twelvepoint$\heartsuit$\kern10pt}}}

 \def\endchange{\InLMargin{\hbox{\sam\twelvepoint$\spadesuit$\kern10pt}}}

 \def\InLMargin#1{\strut\vadjust{%
     \kern-\strutdepth
     \vtop to \strutdepth{%
         \baselineskip\strutdepth
         \llap{\sam$\smash{\hbox{\EdSpider{#1}}}$}\null}}}

 \def\strutdepth{\dp\strutbox}
 \def\strutheight{\ht\strutbox}

 \def\NoteInRMargin#1{\strut\vadjust{%
     \kern-1.001\strutdepth
     \vtop to \strutdepth{%
       \baselineskip\strutdepth
       \vss\rlap{\ninepoint\unskip\hskip\hsize
         \vtop to 0pt{%
           \hsize=16em\hfuzz=\hsize
           \leftskip=10pt%
           \rightskip=0pt plus 10000pt%
           \baselineskip=9.8pt\lineskip=.2pt%
           \let\\\break
           \noindent\EdSpider{#1}\vss}%
                \kern10pt}\hbox{}}
       }}

 \def\ednote#1{\NoteInRMargin{\tentt #1}}

 \def\cbar{\InLMargin{%
      \dimen0=\strutdepth\advance\dimen0 by \lineskip
      \vrule width 3pt
      height \strutheight depth \dimen0 \kern
      3pt}}

 \def\ccbar{\InLMargin{%
      \dimen0=2\strutdepth\advance\dimen0 by 2\lineskip
      \vrule width 3pt
        height 3\strutheight depth \dimen0 \kern
      3pt}}

 \newinsert\TRMargIns
 \dimen\TRMargIns=\maxdimen

  \def\Ednote#1{\insert\TRMargIns{%
       \vbox to 0pt{\hsize=140pt\hfuzz=\hsize
           \leftskip=6pt%
           \rightskip=0pt plus 10000pt%
           \baselineskip=9.8pt\lineskip=.2pt%
           \let\\\break
           \SetPageRemainder
           \vglue540pt\vglue-\PageRemainder
           \noindent\EdSpider{\tentt #1}\vss}%
       \smallskip}}

 \def\KillEdStuff{\def\ednote##1{}\def\Ednote##1{}%
      \let\change\relax\let\beginchange\relax\let\endchange\relax
       \let\cbar\relax\let\ccbar\relax}


  \topskip=12pt
  \newskip\StdBaselineskip 
  \StdBaselineskip 12pt
  \lineskip=1.1pt
  \lineskiplimit=.8pt
  \widowpenalty=10000 
  \clubpenalty=10000  
  \abovedisplayskip=6pt plus 1pt minus 1pt
  \abovedisplayshortskip=3pt plus 1.5pt
  \belowdisplayskip=6pt plus 1pt minus 1pt
  \belowdisplayshortskip=5pt plus 1pt minus 1pt
  \hfuzz=1.5pt   

  \def\StdPretolerance{100}
  \tolerance=\StdPretolerance

  \newdimen\StdMathsurround
  \StdMathsurround=1.5pt 
  \mathsurround=\StdMathsurround
  \Mas                   

   \def\prose{\relax\hbox{\kern.6\StdMathsurround}}
  
  \def\StdParskip{0pt}    
  \parskip=\StdParskip
  \parindent=0.5cm
 

  \def\Times{ptmr  } 
  \def\TimesI{ptmri  } 
  \def\TimesB{ptmb  }
  \def\TimesBI{ptmbi  }
  \def\HelveticaN{phvrrn }

  =\Times at 10bp
  =\TimesB at 10bp
  \font\tenit=\TimesI at 10bp
  =\TimesBI at 10bp

  \font\tenmrm=cmr10  


    =\Times at 9bp 
    \font\nineit=\TimesI at 9bp 
    =\TimesB at 9bp 
    =\TimesBI at 9bp 

    =\HelveticaN at 9bp 


  =\Times at 12bp
  \font\twelveit=\TimesI at 12bp
  =\TimesB at 12bp


  \font\titleit=\TimesI at 14.4bp
  =\TimesB at 14.4bp

 \SetAuthorHead{AuthorHead} 
 \SetTitleHead{TitleHead}  


  \def\lBr{\raise.125ex\hbox{[\kern.1125ex}}
  \def\rBr{\raise.125ex\hbox{\kern.1125ex]}}

 \setbox\footbox=\hbox{\Smallfonts 2)~}



  \bgroup
  \catcode`\@=11 
  \gdef\itSpacing{%
     \xspaceskip=.31em plus.1em minus.05em \sfcode `f=2001
     \itWarning@\let\itWarning@\itWarning@@}
  \gdef\itSpacingOff{%
     \xspaceskip=0pt \sfcode `f=1000
     \let\itWarning@\relax}
   \global\let\itWarning@\relax
  \gdef\itWarning@@{\errmessage{%
  Special italic spacing already in force
  (you have probably omitted an ``endth'').
  See itSpacing macro in osuPSfnt.sty
         }}
  \egroup

 \fontdimen1\titlebf=0.0pt
 \fontdimen2\titlebf=3.6135pt
 \fontdimen3\titlebf=2.8908pt
 \fontdimen4\titlebf=1.44539pt
 \fontdimen5\titlebf=6.64882pt
 \fontdimen6\titlebf=14.45398pt
 \fontdimen7\titlebf=1.60439pt

 \fontdimen1\tenbi=0.26794pt
 \fontdimen2\tenbi=2.50937pt
 \fontdimen3\tenbi=2.00749pt
 \fontdimen4\tenbi=1.00374pt
 \fontdimen5\tenbi=4.59717pt
 \fontdimen6\tenbi=10.03749pt
 \fontdimen7\tenbi=1.11415pt

 \fontdimen1\twelverm=0.0pt
 \fontdimen2\twelverm=3.01125pt
 \fontdimen3\twelverm=2.409pt
 \fontdimen4\twelverm=1.2045pt
 \fontdimen5\twelverm=5.39615pt
 \fontdimen6\twelverm=12.045pt
 \fontdimen7\twelverm=1.33699pt

 \fontdimen1\twelveit=0.27731pt
 \fontdimen2\twelveit=3.01125pt
 \fontdimen3\twelveit=2.409pt
 \fontdimen4\twelveit=1.2045pt
 \fontdimen5\twelveit=5.37207pt
 \fontdimen6\twelveit=12.045pt
 \fontdimen7\twelveit=1.33699pt

 \fontdimen1\twelvebf=0.0pt
 \fontdimen2\twelvebf=3.01125pt
 \fontdimen3\twelvebf=2.409pt
 \fontdimen4\twelvebf=1.2045pt
 \fontdimen5\twelvebf=5.5407pt
 \fontdimen6\twelvebf=12.045pt
 \fontdimen7\twelvebf=1.33699pt

 \fontdimen1\tenrm=0.0pt
 \fontdimen2\tenrm=2.50937pt
 \fontdimen3\tenrm=2.00749pt
 \fontdimen4\tenrm=1.00374pt
 \fontdimen5\tenrm=4.49678pt
 \fontdimen6\tenrm=10.03749pt
 \fontdimen7\tenrm=1.11415pt

 \fontdimen1\tenit=0.27731pt
 \fontdimen2\tenit=2.50937pt
 \fontdimen3\tenit=2.00749pt
 \fontdimen4\tenit=1.00374pt
 \fontdimen5\tenit=4.47672pt
 \fontdimen6\tenit=10.03749pt
 \fontdimen7\tenit=1.11415pt

 \fontdimen1\tenbf=0.0pt
 \fontdimen2\tenbf=2.50937pt
 \fontdimen3\tenbf=2.00749pt
 \fontdimen4\tenbf=1.00374pt
 \fontdimen5\tenbf=4.61723pt
 \fontdimen6\tenbf=10.03749pt
 \fontdimen7\tenbf=1.11415pt

 \fontdimen1\ninerm=0.0pt
 \fontdimen2\ninerm=2.25842pt
 \fontdimen3\ninerm=1.80673pt
 \fontdimen4\ninerm=0.90337pt
 \fontdimen5\ninerm=4.0471pt
 \fontdimen6\ninerm=9.03374pt
 \fontdimen7\ninerm=1.00273pt

 \fontdimen1\nineit=0.27731pt
 \fontdimen2\nineit=2.25842pt
 \fontdimen3\nineit=1.80673pt
 \fontdimen4\nineit=0.90337pt
 \fontdimen5\nineit=4.02904pt
 \fontdimen6\nineit=9.03374pt
 \fontdimen7\nineit=1.00273pt

 \fontdimen1\ninebf=0.0pt
 \fontdimen2\ninebf=2.25842pt
 \fontdimen3\ninebf=1.80673pt
 \fontdimen4\ninebf=0.90337pt
 \fontdimen5\ninebf=4.15552pt
 \fontdimen6\ninebf=9.03374pt
 \fontdimen7\ninebf=1.00273pt


 \newcount\MaxSpaceFactor
 \MaxSpaceFactor=3000 

 \def\ItemStyle{\rm}
 \def\NrStyle{\rm}
 \def\ItemItemStyle{\rm}

 \MaxItemTag{(iii)}
 \MaxItemItemTag{(iii)}
 \MaxNrTag{(2)}
 \MaxFootTag{2)}
 \def\ReferenceHang{30pt}

 \catcode`\@=\active


\loadbold

=\Times  
=\Times scaled750
=\Times scaled650
\font\rms=\Times scaled 920 

=\TimesBI scaled 860
=\TimesI scaled 860

\textfont0=\rrm  
\scriptfont0=\erm 
\scriptscriptfont0=\srm

\def\Augment#1#2{%
    \toks0\expandafter{#1}\toks2{#2}%
    \edef#1{\the\toks0\the\toks2}}

 \font\twelverma=\Times  scaled 1200
 \font\tenrma=\Times  scaled 1000
 \font\ninerma=\Times scaled 920
 =\Times scaled 840
 \font\sevenrma=\Times scaled 760
 =\Times scaled 680
 \font\fiverma=\Times scaled 600

 \Augment\tenpoint{%
  \textfont0=\tenrma  \scriptfont0=\sevenrma  
  \scriptscriptfont0=\fiverma  }

 \Augment\ninepoint{%
  \textfont0=\ninerma  \scriptfont0=\sevenrma 
  \scriptscriptfont0=\fiverma}

 \Augment\twelvepoint{%
  \textfont0=\twelverma  \scriptfont0=\ninerma  
  \scriptscriptfont0=\sevenrma}

\mathsurround=1pt
\hsize=13.45truecm
\vsize=19.5truecm
\hoffset=1.25truecm
\voffset=2truecm
\advance\baselineskip by 2pt

\predefine\til{\~}
\def\~#1{\relax\ifmmode\widetilde{#1}\else\til{#1}\fi}

\redefine \ge{\geqslant}
\define \wt#1{\mathaccent"0365{#1}}
\define \wh#1{\mathaccent"0362{#1}}

\define \iss{\,\Mathaccent{\raise -.8 ex\hbox{$\widetilde{}$\kern.1em}}\rightarrow\,}

\define \prlim{{\varprojlim}\vphantom{i}\,}

\define \ur{\mathop{\fam0 ur}}

\define \dimm{\operatorname{\fam0 dim\,}}

\define \id{\operatorname{\fam0 id\,}}

\define \ab{\mathop{\fam0 ab}}

\define \sep{\mathop{\fam0 sep}}

\define \lln{\operatorname{\fam0 log\,}}

\define \res{\operatorname{\fam0 res}}

\define \Gal{\mathop{\fam0 Gal}}
\define \Hom{\operatorname{\fam0 Hom}}

\define \Rep{\mathop{\fam0 Rep}}

\define \et{\text{\erm \'et}}

\Mas
\HideEdStuff
\rm 
 

\def\issn{{\nineit ISSN 1464-8997 (on line) 1464-8989 (printed)}}

\def\gtp{{\nineit Published 10 December 2000: \ \copyright\ Geometry \& 
Topology Publications}}

\def\gtv3{{\nineit Geometry \& Topology Monographs, Volume 3 (2000) --
Invitation to higher local fields}}


\def\lione
{{\rms Geometry \& Topology Monographs}}

\def \litwo{{\rms Volume 3: Invitation to higher local fields
}} 

\def\tinfo #1.#2.#3-#4
{{
\noindent  {\lione} \hfill 
\par 
\vskip-1.5pt
\noindent {\litwo} \hfill
\par 
\vskip-1,5pt
\noindent {\rms Part #1, section #2, pages #3--#4} \hfill
\vskip24pt 
}}

\def\tinfos #1.#2.#3-#4
{{
\noindent  {\lione} \hfill 
\par 
\vskip-1.5pt
\noindent {\litwo} \hfill
\par 
\vskip-1.5pt
\noindent {\rms Pages #3--#4} \hfill
\vskip24pt 
}}

\def\tinfoi #1
{{
\noindent  {\lione} \hfill 
\par 
\vskip-1.5pt
\noindent {\litwo} \hfill
\par 
\vskip-1.5pt
\noindent {\rms Pages iii--xi: Introduction and contents} \hfill
\vskip26pt 
}}


  \def\titlepagehead{\hfil}

  \newif\iftitlepage\titlepagefalse
  \newif\ifblankpage\blankpagefalse
  \def\makeheadline{
     \ifblankpage{}\else%
     \iftitlepage
\vbox{\line{\vbox to 8.5pt{}
\ninerm
\copy\HLinebox \hfill
\hglue5mm\ninebf\folio 
\titlepagehead}}%
      \else
\vbox{\ifodd\pageno\rightheadline\else\leftheadline\fi}%
      \fi\vskip 12pt\fi}%
     \def\rightheadline{\line{\vbox to 8.5pt{}%
      \ninerm
\copy\TitleBox \hfill
\hglue5mm\ninebf\folio}}%
     \def\leftheadline{\line{\vbox to 8.5pt{}%
        \unskip\ninerm\unskip\ninebf\folio\hglue5mm
 \hfill \copy\AuthorBox
}}

 \footline={\ifblankpage{}\else
\iftitlepage\ninepoint\sam\hfill
\line{\vbox to 8.5pt{}
\copy\TFLinebox
\hfill
\hglue5mm 
}
            \else
\ninepoint\sam\hfill
\line{\vbox to 8.5pt{}
\copy\FLinebox
\hfill 
\hglue5mm
}
\hfil\fi\global\titlepagefalse\fi}

\def\blankpage{{\blankpagetrue\noindent\hbox to 10pt{\hss}\vfill
\pagebreak}}

\tenpoint\rm 
 

\pageno=263

\tinfo II.6.263-272

\SetTFLinebox{\gtp }
\SetFLinebox{\gtv3 }
\SetHLinebox{\issn}

\H 6. $\Phi$-$\Gamma$-modules and Galois cohomology

Laurent Herr

\SetAuthorHead{L. Herr}
\SetTitleHead{Part II. Section  6. $\Phi$-$\Gamma$-modules and Galois cohomology \qquad\qquad}

\HH 6.0. Introduction

Let $G$ be a profinite group and $p$ a prime number.


\df Definition

A finitely generated $\Bbb Z_p$-module $V$
endowed with a continuous $G$-action is called a $\Bbb Z_p$-adic
representation of $G$. Such representations form a category denoted by 
$\Rep_{\,\Bbb Z_p}(G)$; its subcategory $\Rep_{\,\Bbb F_p}(G)$
(respectively $\Rep_{p\text{\erm -tor}}(G)$) of mod $p$ representations (respectively
$p$-torsion representations) consists of the $V$ annihilated by $p$ 
(respectively a power of $p$).
\enddf


\rk Problem

To calculate in a simple explicit way the cohomology groups
$H^i(G,V)$ of the representation $V$.
\endrk

A method to solve it for $G=G_K$ ($K$ is a local field) 
is to use Fontaine's theory of $\Phi$-$\Gamma$-modules and   
pass to a simpler Galois representation, paying the price of enlarging
$\Bbb Z_p$ to the ring of integers of a two-dimensional local field.
In doing this we have to replace linear with semi-linear actions.


In this paper we give an overview of the applications of 
such techniques in different situations. We begin with a simple example. 


\HH 6.1. 
The case of a field of positive characteristic

Let $E$ be a field of characteristic $p$, $G=G_E$
and $\sigma\colon E^{\sep}\to E^{\sep}$, $\sigma(x)=x^p$
the absolute Frobenius map.


\df Definition

For $V\in \Rep_{\Bbb F_p}(G_E)$ put
$D(V):=(E^{\sep}\otimes_{\Bbb F_p}V)^{G_E}$; 
$\sigma$ acts on $D(V)$ by acting on $E^{\sep}$.
\enddf


\rk Properties

\Roster
\Item{(1)} $\dimm_E D(V)=\dimm_{\,\Bbb F_p} V$;

\Item{(2)} the ``Frobenius'' map $\varphi\colon D(V)\to D(V)$ induced by
$\sigma\otimes \id_V$ satisfies: 

\Item{}
\ \ a) $\varphi(\lambda x)=\sigma(\lambda)\varphi(x)$
for all $\lambda\in E$, $x\in D(V)$ (so $\varphi$ is $\sigma$-semilinear);

\Item{}\ \ b) $\varphi(D(V))$ generates $D(V)$ as an $E$-vector space.
\endRoster
\endrk


\df Definition

A finite dimensional vector space $M$ over $E$ is called
an {\it \'etale $\Phi$-module} over $E$ if there is
a $\sigma$-semilinear
 map $\varphi\colon M\to M$ such that
$\varphi(M)$ generates $M$ as an $E$-vector space.
\enddf


\noindent \'Etale $\Phi$-modules form an abelian category
$\Phi M_E^{\et}$ (the morphisms are the linear maps
commuting with the Frobenius $\varphi$).


\th Theorem 1 {{\rm (Fontaine, [F])}}

The functor $V\to D(V)$ is an equivalence of the categories
$\Rep_{\,\Bbb F_p}(G_E)$ and $\Phi M_E^{\et}$.
\endth


We see immediately that $H^0(G_E,V)=V^{G_E}\simeq D(V)^{\varphi}$.

So in order to obtain an explicit description of the Galois
cohomology of mod $p$ representations of $G_E$, we should try
to derive in a simple manner 
the functor associating to an \'etale $\Phi$-module
the group of points fixed under $\varphi$. This is indeed
a much simpler problem because there is only one operator
acting.

For $(M,\varphi)\in \Phi M_E^{\et}$ define the following
complex of abelian groups:
$$C_1(M):\qquad 0@>>> M @>\varphi-1>> M@>>> 0$$
{{\rm(}}$M$ stands at degree 0 and 1{{\rm)}}.

This is a functorial construction, so by taking the cohomology of the
complex, we obtain
a cohomological functor $({\Cal H}^i:=H^i\circ C_1)_{i\in\Bbb N}$
from $\Phi M_E^{\et}$ to the category of abelian groups.


\th Theorem 2

The cohomological functor $({\Cal H}^i\circ D)_{i\in\Bbb N}$
can be identified  with the Galois cohomology functor 
$(H^i(G_E,.\,))_{i\in\Bbb N}$ for the category $\Rep_{\,\Bbb F_p}(G_E)$.
So, if $M=D(V)$ then ${\Cal H}^i(M)$ provides a simple explicit description
of $H^i(G_E,V)$.

\endth


\pf Proof of Theorem 2

We need to check that the cohomological functor
$({\Cal H}^i)_{i\in \Bbb N}$ is universal; therefore
it suffices to verify that for every $i\ge 1$
the functor ${\Cal H}^i$ is effaceable: this means that
for every $(M,\varphi_M) \in \Phi M_E^{\et}$ and every
$x\in {\Cal H}^i(M)$ there exists
an embedding $u$ of $(M,\varphi_M)$ in $(N,\varphi_N)\in 
\Phi M_E^{\et}$ such that
${\Cal H}^i(u)(x)$ is zero in ${\Cal H}^i(N)$.
But this is easy: it is trivial for $i\ge 2$; for $i=1$ 
choose an element $m$ belonging to the class $x\in M/(\varphi-1)(M)$,
put $N:=M\oplus Et$ and extend $\varphi_M$ to 
the $\sigma$-semi-linear map $\varphi_N$ 
determined by $\varphi_N (t):=t+m$.
\qed\endpf


\HH 6.2. $\Phi$-$\Gamma$-modules and $\Bbb Z_p$-adic representations


\df Definition 

Recall that a Cohen ring is an absolutely unramified complete discrete valuation
ring of mixed characteristic $(0,p>0)$, so its  maximal ideal is generated by~$p$.

\enddf


We describe a general formalism, explained by Fontaine in \cite{F},
which lifts the equivalence of categories of Theorem 1
in characteristic $0$ and relates the $\Bbb Z_p$-adic
representations of $G$ to a category of modules over a Cohen ring,
endowed with a ``Frobenius'' map and a group action.

Let $R$ be an algebraically closed complete valuation (of rank 1) field 
of characteristic $p$ and let $H$ be a normal closed subgroup of $G$.
Suppose that $G$ acts continuously on $R$ by ring automorphisms.
Then $F:=R^H$ is a perfect closed subfield of $R$.

For every integer $n\geq 1$, the ring $W_n(R)$ of Witt vectors of
length $n$
is endowed with the product of the topology on $R$ defined by the
valuation
and then $W(R)$ with the inverse limit topology. 
Then the componentwise
action of the group $G$ is continuous and commutes with the natural Frobenius
$\sigma$ on $W(R)$. We also have $W(R)^H=W(F)$.

Let $E$ be a closed subfield of $F$ such that $F$ is the completion of
the $p$-radical closure
of $E$ in $R$. Suppose there exists a Cohen subring 
${\Cal O}_{\Cal E}$ of $W(R)$
with residue field $E$ and which is stable under the actions of $\sigma$ and 
of $G$. 
Denote by $\Cal O_{\widehat {\Cal E}_{\ur}}$ the completion of the 
integral closure of ${\Cal O}_{\Cal E}$ in $W(R)$: 
it is a Cohen ring which
is stable by $\sigma$ and $G$, its residue field is the separable closure
of $E$ in $R$ and 
$(\Cal O_{\widehat {\Cal E}_{\ur}})^H={\Cal O}_{\Cal E}$. 

The natural map from $H$ to $G_E$ is an isomorphism if and only if the action
of $H$ on $R$ induces an isomorphism from $H$ to $G_F$. We suppose that this 
is the case.


\df Definition

Let $\Gamma$ be the quotient group $G/H$. An \'etale $\Phi$-$\Gamma$-module
over $\Cal O_{\Cal E}$ is a finitely generated $\Cal O_{\Cal E}$-module
endowed with a $\sigma$-semi-linear Frobenius map 

\noindent $\varphi\colon M\to M$  
and a 
continuous $\Gamma$-semi-linear action of $\Gamma$ commuting with $\varphi$ 
such that the image of $\varphi$ generates the module $M$.
\enddf


\'Etale $\Phi$-$\Gamma$-modules over $\Cal O_{\Cal E}$ form an abelian 
category 
$\Phi\Gamma M_{{\Cal O}_{\Cal E}}^{\et}$ (the morphisms are
the linear maps commuting with $\varphi$).  There is 
a tensor product of $\Phi$-$\Gamma$-modules, the natural one.
For two objects $M$ and $N$ of $\Phi\Gamma M_{{\Cal O}_{\Cal E}}^{\et}$ the
$\Cal O_{\Cal E}$-module $\Hom_{\Cal O_{\Cal E}}(M,N)$ can be endowed with
an \'etale $\Phi$-$\Gamma$-module structure (see \cite{F}). 

For every $\Bbb Z_p$-adic representation $V$ of $G$, let $D_H(V)$
be the $\Cal O_{\Cal E}$-module 
$(\Cal O_{\widehat {\Cal E}_{\ur}}\otimes_{\Bbb Z_p} V)^H$. It is naturally an
\'etale $\Phi$-$\Gamma$-module,
with $\varphi$ induced by the map $\sigma\otimes \id_V$ and $\Gamma$ acting on
both sides of the tensor product. From Theorem 2 one deduces the 
following fundamental result:


\th Theorem 3 {{\rm (Fontaine, \cite{F})}}

The functor $V\to D_H(V)$ is an equivalence of the categories
$\Rep_{\,\Bbb Z_p}(G)$ and $\Phi\Gamma M_{{\Cal O}_{\Cal E}}^{\et}$.
\endth


\rk Remark

If $E$ is a field of positive characteristic, $\Cal O_{\Cal E}$ is a Cohen
ring with residue field $E$ endowed with a Frobenius $\sigma$, then we can 
easily extend the results of the whole subsection 6.1 to 
$\Bbb Z_p$-adic representations of $G$ by using Theorem 3 
for $G=G_E$ and $H=\{1\}$.

\endrk


\HH 6.3. A brief survey of the theory of the field of norms

For the details we refer to \cite{W}, \cite{FV} or \cite{F}.

Let $K$ be a complete discrete valuation field of characteristic $0$
with perfect residue field $k$  of characteristic $p$. 
Put $G=G_K=\Gal(K^{\sep}/K)$.

Let $\Bbb C$ be the completion of $K^{\sep}$,
denote the extension of the discrete valuation $v_K$ of $K$
to $\Bbb C$ by $v_K$.
Let $R^*=\prlim \Bbb C_n^*$ where $\Bbb C_n=\Bbb C$
and the morphism from $\Bbb C_{n+1}$ to $\Bbb C_n$ is  raising to the 
$p$th power.
Put $R:=R^*\cup\{0\}$ and define
$v_R((x_n))=v_K(x_0)$.
For $(x_n),(y_n)\in R$ define
$$(x_n)+(y_n)=(z_n)\qquad\text{where $z_n=\lim_m\, (x_{n+m}+y_{n+m})^{p^m}$}.$$
Then $R$ is an algebraically closed field   
of characteristic $p$ complete with respect to $v_R$ (cf. \cite{W}). 
Its residue field is isomorphic to the algebraic closure of $k$
and there is a natural continuous action of $G$ on $R$.
(Note that Fontaine denotes this field by $\text{\rrm Fr }R$ in \cite{F}).

Let $L$ be a Galois extension of $K$ in $K^{\sep}$. Recall that
one can always define
the ramification filtration on $\Gal (L/K)$ in the upper numbering.
Roughly speaking, $L/K$ is an arithmetically profinite extension if
one can define the lower ramification subgroups of $G$ so that
the classical relations between the two filtrations  for finite
extensions are preserved. This is in particular possible if
$\Gal (L/K)$ is a $p$-adic Lie group with finite residue field extension.

The field $R$ contains in a natural way the field of norms $N(L/K)$ of every
arithmetically profinite extension $L$ of $K$
and the restriction of $v$ to $N(L/K)$ is a discrete
valuation. The residue field of $N(L/K)$ is isomorphic to that of $L$ and 
$N(L/K)$ is stable under the action of $G$. The construction is 
functorial: if $L'$ is a finite extension of $L$ contained in $K^{\sep}$,
then $L'/K$ is still arithmetically profinite and $N(L'/K)$ is a separable
extension of $N(L/K)$. The direct limit of the fields $N(L'/K)$ where 
$L'$ goes through all the finite extensions of $L$ contained in $K^{\sep}$
is the separable closure $E^{\sep}$ of $E=N(L/K)$. It is stable under
the action of $G$ and the subgroup $G_L$ identifies with $G_E$. The field
$E^{\sep}$ is dense in $R$.

Fontaine described how to lift these constructions in characteristic $0$ 
when $L$ is the cyclotomic $\Bbb Z_p$-extension $K_\infty$ of $K$.
Consider the ring of Witt vectors $W(R)$ endowed with the Frobenius
map $\sigma$ and the natural componentwise action of $G$. Define the topology 
of $W(R)$
as the product of the topology defined by the valuation on
$R$. Then one can construct a Cohen ring $\Cal O_{\widehat {\Cal E}_{\ur}}$
with residue field $E^{\sep}$ ($E=N(L/K)$) such that:

(i) $\Cal O_{\widehat {\Cal E}_{\ur}}$ is stable by $\sigma$ and 
the action of $G$,

(ii) for every finite extension $L$ of $K_\infty$
the ring $(\Cal O_{\widehat {\Cal E}_{\ur}})^{G_L}$ is a Cohen ring with 
residue field $E$.

Denote by ${\Cal O}_{{\Cal E}(K)}$ the ring 
$(\Cal O_{\widehat {\Cal E}_{\ur}})^{G_{K_\infty}}$. It is stable by $\sigma$
and the quotient $\Gamma=G/G_{K_\infty}$ acts continuously on 
${\Cal O}_{{\Cal E}(K)}$ with respect to the induced topology. 
Fix a topological
generator $\gamma$ of $\Gamma$: it is a continuous ring automorphism 
commuting with $\sigma$. The fraction field 
of ${\Cal O}_{{\Cal E}(K)}$ is a two-dimensional standard local field 
(as defined in section 1 of Part~I).
If $\pi$ is a lifting of a prime element of $N(K_\infty/K)$ in
${\Cal O}_{{\Cal E}(K)}$ then the elements of ${\Cal O}_{{\Cal E}(K)}$
are the series $\sum_{i\in\Bbb Z} a_i \pi^i$, where the coefficients
$a_i$ are in $W(k_{K_\infty})$ and converge $p$-adically to $0$
when $i\to -\infty$. 

\HH 6.4. Application  of $\Bbb Z_p$-adic 
representations of $G$ 
\unskip\break \phantom{}\enspace to the Galois cohomology

If we put together Fontaine's construction and the general formalism of
subsection 6.2 we obtain the following important result:


\th Theorem 3' {{\rm (Fontaine, \cite{F})}}

The functor $V\to D(V):=
(\Cal O_{\widehat {\Cal E}_{\ur}}\otimes_{\Bbb Z_p} V)^{G_{K_\infty}}$
defines an equivalence of the categories 
$\Rep_{\,\Bbb Z_p}(G)$ and $\Phi \Gamma M_{{\Cal O}_{{\Cal E}(K)}}^{\et}$.
\endth


Since for every $\Bbb Z_p$-adic representation of $G$
we have $H^0(G,V)=V^G\simeq D(V)^\varphi$, we want now, as in paragraph
6.1, compute explicitly the cohomology of the representation using 
the $\Phi$-$\Gamma$-module associated to $V$.

For every \'etale $\Phi$-$\Gamma$-module $(M,\varphi)$ define the
following complex of abelian groups:
$$C_2(M):\qquad 0@>>> M @>\alpha>> M\oplus M @>\beta>> M@>>> 0$$
where $M$ stands at degree 0 and 2, 
$$\alpha(x)=((\varphi-1)x,(\gamma-1)x),
\quad \beta((y,z))=(\gamma-1)y-(\varphi-1)z.$$
By functoriality, we obtain a cohomological functor 
$({\Cal H}^i:=H^i\circ C_2)_{i\in\Bbb N}$ from 

\noindent $\Phi \Gamma M_{{\Cal O}_{{\Cal E}(K)}}^{\et}$ to the category of abelian
groups.


\th Theorem 4 {{\rm (Herr, \cite{H})}}  

The cohomological functor $({\Cal H}^i\circ D)_{i\in\Bbb N}$
can be identified with the Galois cohomology functor 
$(H^i(G,.\,))_{i\in\Bbb N}$
for the category $\Rep_{p\text{\erm -tor}}(G)$.
So, if $M=D(V)$  then ${\Cal H}^i(M)$ provides a simple explicit description
of $H^i(G,V)$ in the $p$-torsion case.

\endth
 

\pf Idea of the proof of Theorem 4

We have to check that for every $i\ge 1$ the functor $\Cal H^i$ is
effaceable. 
For every $p$-torsion object
$(M,\varphi_M) \in \Phi \Gamma M_{{\Cal O}_{{\Cal E}(K)}}^{\et}$ 
and every
$x\in {\Cal H}^i(M)$ we construct
an explicit embedding $u$ of $(M,\varphi_M)$ in a certain $(N,\varphi_N)\in 
\Phi \Gamma M_{{\Cal O}_{{\Cal E}(K)}}^{\et}$ such that
${\Cal H}^i(u)(x)$ is zero in ${\Cal H}^i(N)$.
For details see \cite{H}. The key point is of topological nature:
we prove, following an idea of Fontaine in \cite{F}, that there exists
an open neighbourhood of $0$ in $M$ on which $(\varphi-1)$ is bijective
and use then the continuity of the action of $\Gamma$. 
\qed\endpf


As an application of theorem 4 we can prove the following result (due to Tate):


\th Theorem 5

Assume that $k_K$ is finite and let $V$ be in $\Rep_{p\text{\erm -tor}}(G)$. 
Without using class field theory the previous theorem implies that
$H^i(G,V)$ are finite, $H^i(G,V)=0$ for $i\ge 3$ and
$$\sum_{i=0}^2 l(H^i(G,V))=-|K\colon \Bbb Q_p|\  l(V),$$
where $l(\ )$ denotes the length over $\Bbb Z_p$.

\endth


See \cite{H}.


\rk Remark

Because the finiteness results imply that the Mittag--Leffler conditions 
are satisfied, it is possible to generalize the explicit
construction of the cohomology and to prove analogous results for
$\Bbb Z_p$ (or $\Bbb Q_p$)-adic representations by passing to the
inverse limits. 
\endrk


\HH 6.5. A new approach to local class field theory

The results of the preceding paragraph allow us to prove without using
class field theory the following:


\th Theorem 6 {{\rm (Tate's local duality)}}

Let $V$ be in $\Rep_{p\text{\erm -tor}}(G)$ and $n\in\Bbb N$ such that 
$p^nV=0$.
Put $V^*(1):=\Hom(V,\mu_{p^n})$. Then there is a canonical isomorphism
from $H^2(G,\mu_{p^n})$ to $\Bbb Z/p^n$ and the cup product 
$$H^i(G,V)\times H^{2-i}(G, V^*(1))@>\cup >> H^2(G,\mu_{p^n})
\simeq \Bbb Z/p^n$$
is a perfect pairing.

\endth


It is well known that a proof of the local duality theorem of Tate
without using class field theory gives a construction of the
reciprocity map. For every $n\ge 1$ we have by duality a
functorial isomorphism between the finite groups $\Hom (G,\Bbb Z/p^n)
=H^1(G,\Bbb Z/p^n)$ and $H^1(G,\mu_{p^n})$ which is isomorphic to
$K^*/(K^*)^{p^n}$ by Kummer theory. Taking the inverse limits
gives us the $p$-part of the reciprocity map, the most difficult part.


\pf Sketch of the proof of Theorem 6 

(\cite{H2}). 


a) Introduction of differentials:

Let us denote by $\Omega^1_c$
the ${\Cal O}_{{\Cal E}(K)}$-module of continuous differential forms
of $\Cal O_{\Cal E}$ over $W(k_{K_\infty})$. 
If $\pi$ is a fixed lifting of a prime element of 
$E(K_\infty/K)$ in $\Cal O_{{\Cal E}(K)}$, then this module is free and
generated by $d\pi$. Define the residue
map from $\Omega^1_c$ to $W(k_{K_\infty})$ by $\res\,\left(\sum_{i\in\Bbb Z}
a_i\pi^id\pi\right):=a_{-1}$; it is independent of the choice of $\pi$.

\smallskip

b) Calculation of some $\Phi$-$\Gamma$-modules:

The ${\Cal O}_{{\Cal E}(K)}$-module $\Omega^1_c$ is endowed with an \'etale
$\Phi$-$\Gamma$-module structure by the following formulas:
for every $\lambda\in {\Cal O}_{{\Cal E}(K)}$ we put: 
$$p\varphi(\lambda d\pi)=
\sigma(\lambda)d(\sigma(\pi))\ \ ,\ \ \gamma(\lambda d\pi)
=\gamma(\lambda) d(\gamma(\pi)).$$

The fundamental fact is that there is a natural isomorphism of 
$\Phi$-$\Gamma$-modules over ${\Cal O}_{{\Cal E}(K)}$ between $D(\mu_{p^n})$
and the reduction $\Omega^1_{c,n}$ of $\Omega^1_c$ modulo $p^n$.

The \'etale $\Phi$-$\Gamma$-module associated to the representation
$V^*(1)$ is 

\noindent $\~M:=\Hom(M,\Omega^1_{c,n})$, where $M=D(V)$. By composing
the residue with the trace we obtain a surjective and continuous map
$\text{\rrm Tr}_n$ from $M$ to $\Bbb Z/p^n $. For every $f\in\~ M$, the map
$\text{\rrm Tr}_n\circ f$ is an element of the group $M^\vee$ of 
continuous group
homomorphisms from $M$ to $\Bbb Z/p^n$. This gives in fact a group
isomorphism from $\~M$ to $M^\vee$ and we can therefore transfer 
the $\Phi$-$\Gamma$-module structure from $\~M$ to $M^\vee$.
But, since $k$ is finite, $M$ is locally compact and $M^\vee$ is in fact the 
Pontryagin dual of $M$. 


\smallskip

c) Pontryagin duality implies local duality:
 
We simply dualize the complex $C_2(M)$ using Pontryagin duality (all arrows
are strict morphisms in the category of topological groups) and obtain
a complex: 
$$C_2(M)^\vee:\qquad
0@>>> M^\vee @>\beta^\vee>> M^\vee\oplus M^\vee @>\alpha^\vee>> M^\vee@>>> 0,$$
where the two $M^\vee$ are in degrees $0$ and $2$.
Since we can construct an explicit quasi-isomorphism between $C_2(M^\vee)$
and $C_2(M)^\vee$, we easily obtain a duality between $\Cal H^i(M)$
and $\Cal H^{2-i}(M^\vee)$ for every $i\in\{0,1,2\}$.


\smallskip

d) The canonical isomorphism from $\Cal H^2(\Omega^1_{c,n})$ to
$\Bbb Z/p^n$:

The map $\text{\rrm Tr}_n$ from $\Omega^1_{c,n}$ to $\Bbb Z/p^n$ factors
through the group $\Cal H^2(\Omega^1_{c,n})$ and this gives an isomorphism.
But it is not canonical! In fact the construction of the complex $C_2(M)$
depends on the choice of $\gamma$. Fortunately, if we take another 
$\gamma$,
we get a quasi-isomorphic complex and if we normalize the map 
$\text{\rrm Tr}_n$
by multiplying it by the unit $-p^{v_p(\lln\chi(\gamma))}/\lln\chi(\gamma)$
of $\Bbb Z_p$, where $\lln$ is the $p$-adic logarithm, $\chi$ the cyclotomic 
character and $v_p=v_{\Bbb Q_p}$, then everything
is compatible with the change of $\gamma$.


e) The duality is given by the cup product:

We can construct explicit formulas for the cup product:  
 $$\Cal H^i(M)\times \Cal H^{2-i}(M^\vee)@>\cup >> \Cal H^2(\Omega^1_{c,n})$$ 
associated with the cohomological functor $(\Cal H^i)_{i\in\Bbb N}$
and we compose them with the preceding normalized isomorphism
from $\Cal H^2(\Omega^1_{c,n})$ to $\Bbb Z/p^n$. Since everything
is explicit, we can compare with the pairing obtained in c) and verify
that it is the same up to a unit of $\Bbb Z_p$.  
\qed\endpf


\rk Remark

Benois, using the previous theorem, deduced an explicit formula of Coleman's 
type  for the Hilbert symbol and proved Perrin-Riou's formula for crystalline
representations (\cite{B}).
\endrk


\HH 6.6. Explicit formulas for the generalized Hilbert symbol 
\unskip\break \phantom{}\enspace on formal groups

Let $K_0$ be the fraction field of the ring $W_0$ of Witt vectors 
with coefficients in a finite field of characteristic $p>2$
and $\Cal F$ a commutative formal group of finite height $h$ defined
over $W_0$.

For every integer $n\ge 1$, denote by $\Cal F[p^n]$ the $p^n$-torsion 
points in $\Cal F(\Cal M_C)$, where $\Cal M_C$ is the maximal ideal
of the completion $C$ of an algebraic closure of $K_0$. 
The group $\Cal F[p^n]$ is isomorphic to $(\Bbb Z/p^n\Bbb Z)^h$.

Let $K$ be a finite extension of $K_0$ contained in $K^{\sep}$
and assume that the points of $\Cal F[p^n]$ are defined over $K$.
We then have a bilinear pairing: 
$$(\ ,\ ]_{\Cal F,n}\colon G_K^{\ab} \times \Cal F(\Cal M_K) \to \Cal F[p^n] $$
(see section~8 of Part~I). 

When the field $K$ contains a primitive $p^n$th root of unity
$\zeta_{p^n}$, Abrashkin gives an explicit description for this
pairing generalizing the classical Br\"uckner--Vostokov formula
for the Hilbert symbol (\cite{A}). In his paper he notices that
the formula makes sense even if $K$ does not contain $\zeta_{p^n}$
and he asks whether it holds without this assumption. In a recent 
unpublished work, Benois  proves that this is true.

Suppose for simplicity that $K$ contains only $\zeta_p$. Abrashkin considers
in his paper the extension $\~K:=K(\pi^{p^{-r}}, r\ge 1)$, where
$\pi$ is a fixed prime element of $K$. It is not a Galois extension of
$K$ but is arithmetically profinite,
so by \cite{W}  one can  consider the field of norms 
for it. 
In order not to loose  information given by the roots of unity of order
a power of $p$, Benois uses the composite
Galois extension $L:=K_\infty \~K/K$ which is arithmetically 
profinite. There are several problems with the field of norms
$N(L/K)$, especially it is not clear that one can lift it in characteristic
$0$ with its Galois action. So, Benois simply considers the completion $F$ 
of the $p$-radical closure of $E=N(L/K)$ and its separable closure 
$F^{\sep}$ in $R$.
If we apply what was explained in subsection 6.2 for $\Gamma=\Gal(L/K)$, 
we get: 


\th Theorem 7 

The functor $V\to D(V):=
(W(F^{\sep})\otimes_{\Bbb Z_p} V)^{G_L}$
defines an equivalence of the categories  
$\Rep_{\,\Bbb Z_p}(G)$ and $\Phi \Gamma M_{W(F)}^{\et}$.
\endth


Choose a topological generator $\gamma'$ of $\Gal(L/K_\infty)$ and
lift $\gamma$ to an element of $\Gal (L/\~K)$. Then $\Gamma$
is topologically generated by $\gamma$ and $\gamma'$, with the relation
$\gamma\gamma'=(\gamma')^a\gamma$, where $a=\chi(\gamma)$
($\chi$ is the cyclotomic character).
For $(M,\varphi)\in \Phi \Gamma M_{W(F)}^{\et}$ the continuous action of
$\Gal(L/K_\infty)$ on $M$ makes it a module over the Iwasawa algebra
$\Bbb Z_p[[\gamma'-1]]$. So we can define the following complex
of abelian groups:
$$C_3(M):\qquad\qquad\qquad 0@>>> M_0@>\alpha\mapsto A_0\alpha >> M_1@>\alpha
\mapsto 
A_1\alpha>> M_2@>\alpha\mapsto A_2\alpha>>M_3 @>>> 0$$
where $M_0$ is in degree $0$,
$M_0=M_3=M$, $M_1=M_2=M^3$,
$$
A_0=\pmatrix \varphi-1\\ \gamma-1\\ \gamma'-1\endpmatrix , 
A_1=\pmatrix \gamma-1 & 1-\varphi  & 0\\
\gamma'-1 & 0 & 1-\varphi \\
0 & {\gamma'}^a -1 & \delta-\gamma
\endpmatrix, 
 A_2=((\gamma')^a-1\ \ \delta-\gamma\ \ \varphi-1)
$$
and
$\delta=((\gamma')^a-1)(\gamma'-1)^{-1} \in \Bbb Z_p[[\gamma'-1]]$.

As usual, by taking the cohomology of this complex, one defines a cohomological
functor $(\Cal H^i)_{i\in\Bbb N}$ from $\Phi \Gamma M_{W(F)}^{\et}$ in the
category of abelian groups. Benois proves only that the cohomology of
a $p$-torsion representation $V$ of $G$ injects in the groups
$\Cal H^i(D(V))$ which is enough to get the explicit formula. 
But in fact a stronger statement  is true:


\th Theorem 8

The cohomological functor $({\Cal H}^i\circ D)_{i\in\Bbb N}$
can be identified  with the Galois cohomology functor 
$(H^i(G,.\,))_{i\in\Bbb N}$
\, for the category $\Rep_{p\text{\erm -tor}}(G)$.

\endth
 

\pf Idea of the proof

Use the same method as in the proof of Theorem 4. It is only
more technically complicated because of the structure of $\Gamma$.
\qed
\endpf


Finally, one can explicitly construct the cup products in terms of the groups
$\Cal H^i$ and, as in \cite{B},  Benois uses them to calculate the 
Hilbert symbol.  


\rk Remark

Analogous constructions (equivalence of category, explicit construction
of the cohomology by a complex) seem to work for higher dimensional
local fields. In particular, in the two-dimensional case, the formalism
is similar to that of this paragraph; the group $\Gamma$ 
acting on the $\Phi$-$\Gamma$-modules has the same
structure as here and thus the complex is of the same form.
This work is still in progress.
\endrk

\Bib References 


\rf{A} V. A. Abrashkin, Explicit formulae for the Hilbert symbol of
a formal group over the Witt vectors, Izv. RAN Math. 61(1997),  463--515.

\rf{B} D. Benois, On Iwasawa theory of crystalline representations, Duke
Math. J. 104 (2) (2000),  211--267.

\rf{FV} I. Fesenko and  S. Vostokov, Local Fields and Their Extensions,   
Trans. of Math. Monographs, 121, A.M.S., 1993.


\rf{F} J.-M. Fontaine, Repr\'esentations $p$-adiques des corps locaux,
The Grothendieck Festschrift 2, Birkh\"auser, 1994, 59--111.


\rf{H1} L. Herr, Sur la cohomologie galoisienne des corps $p$-adiques,
Bull. de la Soc. Math. de France, 126(1998), 563--600.

\rf{H2} L. Herr, Une approche nouvelle de la dualit\'e locale de Tate, to
appear in  Math. Annalen.

\rf{W} J.-P. Wintenberger, Le corps des normes de certaines extensions
infinies de corps locaux; applications, Ann. Sci. E.N.S. 16(1983), 
59--89.

\endBib

\Coordinates 

Universit\'e Bordeaux 1, Laboratoire A2X, 

351 Cours de la Lib\'eration,
33405 Talence, France

E-mail: herr\@math.u-bordeaux.fr
\endCoordinates

\vfill
\pagebreak

\end